\input francais.sty
\selectlanguage{francais}
\frenchspacing

\input BOXEDEPS.TEX
\input BOXEDEPS.CFG
\HideDisplacementBoxes

\mag=1250
\hsize 16truecm  
\hoffset 0.05truecm
\vsize 22.2truecm
\voffset -0.04truecm
\baselineskip=14pt plus 1pt
\parskip=5pt plus 1pt
\def\cqfd{\unskip\kern 6pt{\penalty 500
\raise -2pt\hbox{\vrule\vbox to 10pt{\hrule width 4pt
\vfill\hrule}}\vrule}\par}

\def\decale#1{\par\noindent\hskip 3em\llap{#1\enspace}\ignorespaces}

{\nopagenumbers
\centerline{\bf N\oe uds non concordants \`a un ${\bf C}$-bord.}
\centerline{M. BOILEAU et L. RUDOLPH\footnote{$^1$}{Partially 
supported by NSF grant DMS-9504832.}}
\vskip2cm

\noindent
{\it RESUME}\par

Un entrelacs orient\'e $L$ dans $S^3$ est un  
${\bf C}$-bord s'il borde dans la boule $B^4$ un\break%
morceau de courbe alg\'ebrique. En utilisant la preuve 
r\'ecente de Kronheimer et Mrowka de la Conjecture de Thom, 
on construit beaucoup de nceuds orient\'es qui ne
sont pas concordants a un ${\bf C}$-bord. Dans le cas 
ou $S^3$ est le bord strictement pseudo\-convexe de $B^4$, 
en utilisant Ie polynome \`a deux variables HOMFLY, on donne une
obstruction pour un entrelacs \`a etre un ${\bf C}$-bord.
\bigskip

\noindent
{\it ABSTRACT}\par

An oriented knot $L$ in $S^3$ is a ${\bf C}$-boundary if it 
bounds a piece of algebraic curve in the ball $B^4$. Using 
Kronheimer and Mrowka's recent proof of the Thom Conjecture, 
we construct many oriented knots which are not concordant to 
a ${\bf C}$-boundary. We use the two-variable HOMFLY polynomial 
to give an obstruction to a knot's being a ${\bf C}$-boundary 
in a strictly pseudoconvex $S^3$.
\vfill
\eject
}
\pageno=1

\centerline {\bf N\oe uds non concordants \`a un $\bf C$-bord.} 
\centerline {M. BOILEAU et L. RUDOLPH}
\bigskip
\bigskip
\bigskip 
\noindent \S\kern.15em 0 - {\bf Introduction.} 
\medskip
Un entrelacs est une r\'eunion finie de cercles 
disjoints plong\'es de fa\c con lisse dans la 
sph\`ere $S^3$. Un n\oe ud est un entrelacs connexe. \par

En identifiant $S^3$ avec le bord d'une sph\`ere ronde 
dans ${\bf C}^2, \{(z, \omega) \in {\bf C}^2 \mid 
|z|^2 + |\omega|^2 = R^2\}$ on peut construire beaucoup
d'entrelacs int\'eressants comme intersection transverse 
de $S^3$ avec une courbe alg\'ebrique plane, par 
exemple : les entrelacs alg\'ebriques des singularit\'es
isol\'ees de courbes alg\'ebriques planes lorsque 
$R << 1$ $([{\rm EN}], [{\rm Le}], [{\rm Mi}]),$
les entrelacs \`a l'infini de courbes alg\'ebriques 
planes lorsque $R >> 1$ $([{\rm Ne}], [{\rm NR}], 
[{\rm R_1}])$, mais aussi beaucoup d'autres entrelacs
comme les
entrelacs admettant une pr\'esentation sous forme de tresses ferm\'ees
quasipositives
$([{\rm R_2}], [{\rm R_3}], [{\rm R_4}], [{\rm R_5}], [{\rm R_6}], [{\rm
R_7}], [{\rm
R_8}])$. \par

Ceci conduit \`a donner la d\'efinition suivante :
\bigskip

\noindent {\bf 0.1.~~D\'efinition.}~~{\it On dit 
qu'un entrelacs orient\'e $L$ dans $S^3$
est un ${\bf C}$-bord s'il existe une boule 
lisse $B^4 \subset {\bf C}^2$ et une courbe
alg\'ebrique plane \'eventuellement singuli\`ere 
$V \subset {\bf C}^2$ qui rencontre
transversallement $\partial B^4$ 
et telle que les paires $(S^3, L)$ 
et~$(\partial B^4,\partial B^4 \cap V)$  
soient diff\'eomorphes par un diff\'eomorphisme
pr\'eservant les ori\-entations ; la 
paire $(\partial B^4, \partial B^4 \cap V)$
\'etant naturellement orient\'ee par l'orientation 
complexe de $B^4$ et de $V$.} \par

Dans $([R_4])$ le second auteur a montr\'e que l'on peut r\'ealiser toute
matrice de Seifert d'un entrelacs orient\'e, comme matrice de Seifert d'un
entrelacs quasipositif, en particulier d'un ${\bf C}$-bord. Il s'en suit
qu'aucun
invariant calculable \`a partir d'une matrice de Seifert ne peut dire si un
entrelacs
orient\'e est un ${\bf C}$-bord, ni m\^eme s'il est concordant \`a un ${\bf
C}$-bord.
\par

Le but de cette note est d'utiliser la preuve r\'ecente par Krohneiner et
Mrowka
$[KM_{1,2}]$ (cf. aussi $[{\rm R_{12}}]$, $[{\rm R_{13}}]$) de la
conjecture de Thom
locale (ou bien, leur preuve plus r\'ecente - et peut-\^etre plus simple -
de la
conjecture de Thom globale $([{\rm KM_3}]))$ pour construire beaucoup de
n\oe uds
orient\'es qui ne sont pas concordants \`a un n\oe ud ${\bf C}$-bord.
\vfill\eject

\noindent {\bf 0.2. Remarques.} \par
{\bf 1)} Dans la d\'efinition d'un ${\bf C}$-bord, on n'impose pas \`a la
boule lisse
$B^4 \subset {\bf C}^2$ d'avoir un bord strictement pseudo-convexe. \par
Dans le cas o\`u le bord de la boule $B^4$ est strictement pseudo-convexe,
on dira
que l'entrelacs orient\'e $L$ est un {\it spc}-${\bf C}$-bord. C'est le cas de
tous les
entrelacs ${\bf C}$-bords connus actuellement, et en particulier des
entrelacs
alg\'ebriques, des entrelacs \`a l'infini des courbes alg\'ebriques planes
et des
entrelacs quasipositifs. \par
Il est important de noter que l'existence d'invariants, permettant de
d\'ecider si un
n\oe ud ou entrelacs orient\'e est un ${\bf C}$-bord, reste un probl\`eme ouvert. Par
contre, le~polyn\^ome \`a $2$ variables HOMFLY permet de donner une
condition
n\'ecessaire pour qu'un entrelacs orient\'e soit un 
{\it spc}-${\bf C}$-bord (cf.
\S\kern.15em 3). Cette condition \'etait connue pour~les entrelacs
quasipositifs
$[{\rm R_7}], [{\rm R_9}], [{\rm R_{10}}]$. \par  

{\bf 2)} Lorsque $\partial B^4$ est strictement pseudo-convexe, on dit
qu'un
plongement lisse \ $\gamma$ \ d'une famille disjointe de cercles dans
$\partial B^4$
est polynomialement convexe si son enveloppe convexe polynomiale dans ${\bf
C}^2$ est
r\'eduite \`a $\gamma$. Avec cette \hbox{terminologie,} 
d'apr\`es Stolzenberg [St], un {\it spc}-${\bf C}$-bord est 
un entrelacs qui admet un plongement lisse dans $\partial
B^4$ qui
n'est pas polynomialement convexe. Un des buts de cette note est~donc de
construire
des n\oe uds dans $\partial B^4$ pour lesquels {\bf tout plongement} est
\hbox{polynomialement convexe}. En g\'en\'eral un plongement d'un n\oe ud dans
$\partial B^4$ est\break%
seulement g\'en\'eriquement polynomialement convexe (cf. [Sch]). 

\vskip 1truecm

\noindent \S\kern.15em 1 - {\bf C-bord et grand genre de Murasugi.} 
\bigskip

Le but de ce paragraphe est d'\'etudier le comportement du grand genre de
Murasugi
(cf. [BW]) des ${\bf C}$-bord sous l'op\'eration de somme connexe. \par

Nous commen\c cons par rappeler la d\'efinition d'une concordance entre
deux entrelacs
orient\'es, ainsi que celle du grand genre de Murasugi.
\medskip 

\noindent {\bf 1.1. D\'efinition.} {\it Deux entrelacs orient\'es $L_0$ et
$L_1$ sont
dits concordants s'il existe un plongement lisse $\Phi :
\displaystyle{\coprod^r_{i=1}} S^1_i \times [0,1] \to S^3 \times [0, 1]$
tel que : 

\item {\bf i)} $\Phi^{-1}(S^3 \times \{0\}) =
\displaystyle{\coprod^r_{i=1}} S^1_i
\times \{0\} = L_0$ \par 
\eject

$\Phi^{-1}(S^3 \times \{1\}) = \displaystyle{\coprod^r_{i=1}} S^1_i \times
\{1\} =
L_1$

\item {\bf ii)} le bord orient\'e de la paire $(S^3 \times [0, 1], Im
\Phi)$ est
$(S^3 \times \{0\}, L_0) \amalg$ $(S^3 \times \{1\}, - L_1)$.}

On sait calculer \`a partir d'une forme de Seifert d'un entrelacs orient\'e
des
invariants de concordance, comme la signature. Dans cette note nous
utiliserons un
invariant plus g\'eom\'etrique, mais aussi plus difficile \`a calculer, qui
est le
grand genre de Murasugi (ou grand genre dans la boule $B^4)$ d'un entrelacs
orient\'e.
\par

Etant donn\'e un entrelacs orient\'e $L$ dans $S^3$ vu comme le bord de la
boule
$B^4$, il existe des surfaces compactes orient\'ees $\Sigma^2$ proprement
plong\'ees
de fa\c con lisse dans $B^4$ et dont le bord orient\'e est $L$. Une telle
surface est
appel\'ee surface de Murasugi lisse de l'entrelacs orient\'ee $L \ ([{\rm
Mu}]$, cf.
aussi [BW]). \par

Pour une surface \`a bord $\Sigma^2$, on appelle grand genre de $\Sigma^2$,
et on note
$G(\Sigma^2)$, le genre de la surface ferm\'ee obtenue en collant un disque
trou\'e le
long des composantes du bord de $\Sigma^2$. On a alors la relation suivante
avec
la caract\'eristique d'Euler de $\Sigma^2 : \chi(\Sigma^2) = r -
2G(\Sigma^2)$,
o\`u $r$ est le nombre de composantes connexes de bord de $\Sigma^2$.
\medskip

\noindent {\bf 1.2. D\'efinition.} {\it Soit $L$ un entrelacs orient\'e
dans $S^3$, on
appelle grand genre de Murasugi de $L$, et on note $M(L)$, le plus petit
des
grands genres des surfaces de Murasugi lisses pour $L$. (cf. [BW]).} \par

Clairement, le grand genre de Murasugi d'un entrelacs orient\'e est un
invariant de la
classe de concordance lisse de cet entrelacs. Dans le cas d'un n\oe ud, il
s'agit du~genre de Murasugi classique (cf. [Mu]). \par

R\'ecemment, P. Kronheimer et T. Mrowka ([KM$_{1,2}$], cf. aussi [KM$_3$])
ont
d\'emon-\break tr\'e le r\'esultat suivant : 
\medskip

\noindent {\bf 1.3. Th\'eor\`eme $([KM_{1,2}])$.} {\it Soient $V$ une
courbe
alg\'ebrique lisse dans ${\bf C}^2$ transverse au bord d'une boule $B^4$
plong\'ee de
fa\c con lisse dans ${\bf C}^2$ et $L = \partial B^4 \cap V$ un ${\bf
C}$-bord dans
$\partial B^4$. Alors $M(L) = G(V \cap B^4)$.} \qquad \qquad \qquad \cqfd  
\medskip

Une application surprenante de ce r\'esultat est la propri\'et\'e suivante qui est
cruciale pour construire des entrelacs qui ne sont pas concordants \`a des
${\bf
C}$-bords. \medskip

\noindent {\bf 1.4.~Proposition.} {\it Soit $L_1 \# L_2$ la somme connexe
de deux
entrelacs orient\'es concordants \`a des  ${\bf C}$-bords. Alors
$M(L_1 \# L_2) \geq M(L_1) + M(L_2) - 1$.}
\medskip

{\bf Preuve.} Puisque l'op\'eration de somme connexe est compatible avec la
relation
de concordance, il suffit de d\'emontrer la proposition dans le cas o\`u
$L_1$ et
$L_2$ sont des ${\bf C}$-bords. Soient $L_1 =
\displaystyle{\bigcup^{r_1}_{i=1}}
K_{i,1}$ et $L_2 = \displaystyle{\bigcup^{r_1}_{j=1}} K_{j,2}$ deux ${\bf
C}$-bords.
Soit $L_1 \# L_2$ la somme~connexe de $L_1$ et $L_2$ telle que $L_1 \# L_2
=
\bigl(K_{1,1} \# K_{1,2}\bigr) \cup  \displaystyle{\bigcup^{r_1}_{i=1}}
K_{i,1} \cup
\displaystyle{\bigcup^{r_2}_{j=2}}  K_{j,2}$. \par

On peut alors construire un cobordisme lisse $W^2$ dans $S^3 \times [0, 1]$
entre
$L_1 \# L_2$ dans $S^3 \times \{0\}$ et l'entrelacs scind\'e $L_1 \amalg
L_2$ dans
$S^3 \times \{1\}$. Ce cobordisme, proprement plong\'e de fa\c con lisse
dans $S^3
\times [0,1]$, est la r\'eunion disjointe d'une surface planaire ayant 3
composantes
de bord $(K_{1,1} \# K_{1,2}) \cup K_{1,1} \cup K_{1,2}$ et d'anneaux
joignant les
autres composantes de $L_1 \# L_2$ aux composantes correspondantes de $L_1
\amalg
L_2$. (cf. Figure  o\`u $r_1 = r_2 = 1)$.

\centerline{\BoxedEPSF{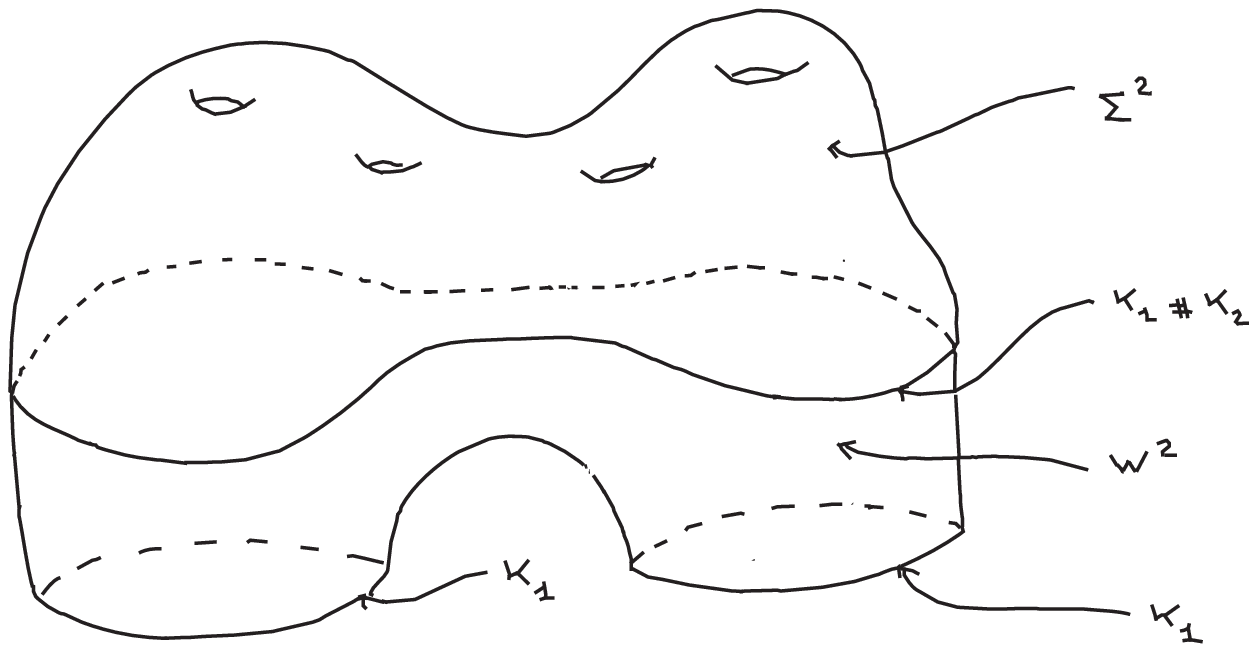}}


Soit $\Sigma^2$ une surface lisse proprement plong\'ee dans $B^4$, de bord
$L_1$ \#
$L_2$ et de grand genre $G(\Sigma^2) = M(L_1 \# L_2)$. Alors, la surface
$\Sigma^2
\displaystyle{\bigcup_{L_1 \# L_2}} W^2$ borde l'entrelacs orient\'e
scind\'e $L_1
\amalg L_2$ dans $B'^4 = B^4  \displaystyle{\bigcup_{\partial
B^4=S^3\times\{0\}}} S^3
\times [0,1]$. Il en d\'ecoule imm\'ediatement l'in\'egalit\'e :
$M(L_1 \# L_2) = G(\Sigma^2) \geq M(L_1 \amalg L_2)-1$. \par

Puisque $L_1$ et $L_2$ sont des ${\bf C}$-bords, l'entrelacs orient\'e
scind\'e $L_1
\amalg L_2$ est un ${\bf C}$-bord, bord de la r\'eunion (disjointe) des
deux
morceaux de courbes alg\'ebriques planes bord\'ees par $L_1$ et $L_2$ dans
$B^4$. Le
r\'esultat de Kronheimer et Mrowka (Theorem. 1.3)  entra\^\i ne alors que
$M(L_1 \amalg
L_2) = M(L_1) + M(L_2)$. D'o\`u l'in\'egalit\'e $M(L_1 \# L_2) \geq M(L_1)
+
M(L_2)-1$. \qquad \qquad \cqfd
\medskip

\noindent {\bf 1.5. Remarque :} 
L'in\'egalit\'e obtenue dans la proposition 1.4 donne une grande
restriction sur les
classes de concordances des ${\bf C}$-bords. En g\'en\'eral, $M(L)-1$ est
loin d'\^etre suradditif par somme connexe. Un exemple frappant (qui sera
exploit\'e
dans la construction des n\oe uds non concordants \`a un ${\bf C}$-bord)
est le fait
que pour un n\oe ud quelconque $K$, $M(K \# K^*) = 0$, o\`u $K^*$ est
obtenu en prenant
l'image de $K$ dans un miroir et en renversant son orientation (cf. 2.1).
En
particulier, on a le corollaire imm\'ediat de la proposition 1.4.
\medskip

\noindent {\bf 1.6. Corollaire.} {\it Soit $K$ un ${\bf C}$-bord, alors $K$
est soit
nul concordant, soit d'ordre infini dans le groupe de concordance lisse des
n\oe uds
orient\'es de $S^3$.}
\medskip

{\bf Preuve.} Soit $K$ un ${\bf C}$-bord. Si $K$ est d'ordre fini $p$ dans
le groupe de
concordance lisse des n\oe uds orient\'es de $S^3$, alors $ 0 = M(K \# K \#
\cdots \#
K) \geq pM(K) - (p-1)  ;$ d'o\`u $M(K) = 0$ et $K$ est nul concordant.
\qquad \qquad
\cqfd 
\medskip

\noindent {\bf 1.7. Remarque.} Dans $\rm{[R_{12}]}$, le second auteur
montre que les
doubles it\'er\'es non tordus positifs, d'une tresse positive sont des
n\oe uds fortement quasipositifs, en particulier des ${\bf C}$-bords. De
m\^eme les
n\oe uds de bretzel $K(p, q, r)$ tels que $p, q, r$ sont premiers,
diff\'erents de
$\pm 1$ et v\'erifient $qr + rp + pq = -1$, sont des ${\bf C}$-bords, (cf.
$\rm{[R_{12}])}$. Ces n\oe uds ont la propri\'et\'e d'avoir un polyn\^ome
d'Alexander
trivial et donc, d'apr\`es M. Freedman [F], sont topologiquement nul
concordants,
c'est-\`a-dire bordent dans $B^4$ un disque topologique localement plat, proprement
plong\'e. \par

Dans $\rm{[R_{12}] [R_{13}]}$, le second auteur montre que ces n\oe uds ne
sont pas nul
con\-cordants~de~fa\c con lisse. Il d\'ecoule du Corollaire 1.6. que tous ces
n\oe uds
sont en fait d'ordre infini dans le groupe de concordance lisse. \par

En fait, il existe un homomorphisme naturel du groupe de concordance lisse
${\bf{\cal K}}_{Diff}$ dans le groupe de concordance topologique localement
plate
${\bf{\cal K}}_{Top}$, dont le noyau contient des \'el\'ements d'ordre
infini
(cf.[CG], [Go]). Le Corollaire 1.6 permet ainsi d'exhiber beaucoup
d'\'el\'ements
d'ordre infini dans ce noyau. La question suivante est \`a notre
connaissance toujours
ouverte (cf.[CG], [Go]). \medskip

{\bf Question 1.} {\it Le noyau de l'homomorphisme naturel $h : {\bf{\cal
K}}_{Diff}
\to {\bf{\cal K}}_{Top}$ \break contient-il un sous-groupe ab\'elien de
rang infini ?}
\par Cette question est li\'ee \`a la question suivante :
\eject

{\bf Question 2.} {\it Deux n\oe uds $K$ et $K'$, fortement quasipositifs
et qui
sont concordants dans la cat\'egorie lisse, sont-ils \'equivalents (i.e.
existe-t-il
un diff\'eomorphisme pr\'eservant l'orientation de $S^3$ et envoyant $K$
sur $K')?$}
\smallskip

\noindent {\bf 1.8. Remarque.} Une r\'eponse positive \`a la question 2
implique une
r\'eponse positive \`a la question 1. \par
Rappelons qu'un n\oe ud ou entrelacs orient\'e $L$ est dit fortement
quasipositif
s'il est le bord d'une surface de Seifert dans $S^3$ obtenue comme surface
tress\'ee
(ou de Markov) \`a partir d'une pr\'esentation en tresse ferm\'ee
quasipositive de
$L$ (cf. ${\rm [R_3], [R_6], [R_8], [R_{11}]).}$ \par  

Une r\'eponse positive \`a la question 2 implique que le noyau de
l'homomorphisme $ h
: {\bf{\cal K}}_{Diff} \to {\bf{\cal K}}_{Top}$ contient le sous-groupe
ab\'elien de
rang infini engendr\'e par tous~les~n\oe uds fortement quasipositifs de
polyn\^ome
d'Alexander trivial. En effet, d'apr\`es le Corollaire 1.6 tous ces n\oe
uds sont
d'ordre infini dans ${\bf{\cal K}}_{Diff}$ (cf. aussi ${\rm[R_{12}],
[R_{13}]).}$ Ils
sont de plus lin\'eairement ind\'ependants puisque la somme connexe de deux
n\oe uds
fortement quasipositifs est un n\oe ud fortement quasipositif. D'autre
part cette
famille de n\oe uds est infini car elle contient tous les doubles positifs
non tordus
des n\oe uds toriques $T_{p,q}$ (cf ${\rm [R_{12}]).}$ \par

Pour un entrelacs $L$ fortement quasipositif, on sait que son miroire $L^*$
n'est
jamais quasipositif. Dans le cas des ${\bf C}$-bords on a :
\smallskip

\noindent {\bf 1.9. Corollaire.} {\it Soit $L$ un ${\bf C}$-bord dans
$S^3$. Si le
miroire $L^*$ de $L$ est un ${\bf C}$-bord, alors $L$ est un entrelacs nul
concordant, c'est-\`a-dire $M(L) = 0$.}
\smallskip

{\bf Preuve.} Si $L$ et $L^*$ sont des ${\bf C}$-bords, alors : \par
\noindent $0 = M(L \# L^*) \geq M(L) + M(L^*) - 1 = 2M(L) - 1 ;$ d'o\`u
$M(L) = 0$. \qquad \qquad
\cqfd 

\vskip .8truecm
 
\noindent \S\kern.15em 2 - {\bf N\oe uds non concordants \`a un n\oe ud
C-bord.}  
\medskip

Dans cette section nous ne consid\'erons que des n\oe uds orient\'es. 
\medskip

\noindent {\bf 2.1. D\'efinition.} {\it Soit $K$ un n\oe ud orient\'e dans
$S^3$, on
note $K^*$ le n\oe ud orient\'e dans $S^3$ obtenu en prenant l'image de $K$
par un
diff\'eomorphisme de $S^3$ renversant l'orientation de $S^3$ et en
renversant son
orientation.} \par
On rappelle l'op\'eration de somme connexe par bande de deux n\oe uds
orient\'es $K_1$
et $K_2$.
\eject

\noindent {\bf 2.2. D\'efinition.} {\it Soient $K_1$ et $K_2$ deux n\oe uds
orient\'es
plong\'es dans $S^3$ de fa\c con qu'il~existe une sph\`ere $S^2$ qui les
s\'epare.
Soient l'intervalle $I = [0,1]$ et $b : I \times I \to S^3$ un plongement
lisse tel
que $b^{-1}(K_1) = \{0\} \times [0,1]$, $b^{-1}(K_2) = \{1\} \times [0,1]$
et les
orientations sur $K_1$ et $K_2$ sont compatibles avec l'orientation de $I
\times I$.
On appelle somme connexe par bande $b$ de $K_1$ et $K_2$, le n\oe ud
orient\'e}
$$K_1 \# _b K_2 = \bigl\{K_1 \amalg K_2 - b\bigl(b^{-1}(K_1 \amalg
K_2)\bigr)\bigr\}
\cup b(I \times I \bigr\}.$$

Dans le cas o\`u l'image de la bande Imb rencontre la sph\`ere $S^2$ qui
s\'epare
$K_1$ de $K_2$ uniquement le long d'un arc de la forme $b(* \times I)$, on
obtient la
somme connexe usuelle $K_1 \# K_2$ de $K_1$ et $K_2$.  Quelle que soit la
bande $b$,
le n\oe ud orient\'e $K_1 \#_b K_2$ est concordant \`a $K_1 \# K_2$. Mais,
en
g\'en\'eral $K_1 \#_b K_2$ est un n\oe ud premier. \par Nous pouvons
maintenant
\'enoncer : \bigskip

\noindent {\bf 2.3. Th\'eor\`eme.} {\it Soit $K$ un n\oe ud orient\'e
quelconque dans
$S^3$. Soit $J$ un n\oe ud orient\'e concordant \`a un $\bf {C}$-bord,  et
tel que
$M(J) > M(K)$. Alors aucune somme connexe par bande $K \#_b J^*$ n'est
concordante \`a
un $\bf {C}$-bord.} \bigskip

{\bf Preuve.} Puisque $K \#_b J^*$ est concordant \`a $K \# J^*$, il suffit
de
d\'emontrer le th\'eor\`eme pour $K \# J^*$. Supposons que $K \# J^*$ est
concordant
\`a un $\bf {C}$-bord. Alors, puisque $J$ est concordant \`a un $\bf
{C}$-bord, la
proposition 1.4 entra\^\i ne l'in\'egalit\'e : \hfill\break $M\bigl( (K \#
J^*) \#
J\bigr)  \geq M(K \# J^*) + M(J) - 1$. \par
D'autre part, comme la somme connexe est associative, $(K \# J^*) \# J = 
\break K \#
(J^* \# J)$. Or $J^* \# J$ est nul concordant, d'o\`u $K \# (J^* \# J)$ est
concordant
\`a $K$ et $M\bigl( (K \# J^*) \# J\bigr) = M(K)$. \par
On obtient donc $M(K) \geq M(J) + M(K \# J^*)-1$, ce qui contredit
l'hypoth\`ese,
puisque $M(K \# J^*) \geq 1$ et $M(J) > M(K)$. \qquad \qquad \cqfd
\bigskip

La construction donn\'ee par le Th\'eor\`eme 2.3 permet de montrer :
\bigskip

\noindent {\bf 2.4.~~Corollaire.}~~{\it Il existe des n\oe uds orient\'es,
premiers,
d'ordre infini dans le groupe de concordance lisse tels que ni $K$, ni
$K^*$ ne sont
concordants \`a un $\bf {C}$-bord.}
\bigskip

{\bf Preuve.} Soient $K_1$ et $K_2$ deux $\bf {C}$-bords tels que $M(K_1) =
M(K_2)$
et \break  $M(K_1 \# K^*_2) \geq 2$.  D'apr\`es $[KL]$ il existe un n\oe
ud $K$ premier concordant \`a $K_1 \# K^*_2$. La d\'emonstration du
Th\'eor\`eme 2.3
montre que $K$ n'est pas concordant \`a un $\bf {C}$-bord. De m\^eme,
puisque $K^*$
est concordant \`a $K^*_1 \# K_2$, la d\'emonstration du Th\'eor\`eme 2.3
implique que
$K^*$ n'est pas concordant \`a un $\bf {C}$-bord.  \qquad \qquad \cqfd
\bigskip

Nous allons construire d'autres exemples de n\oe uds non concordants \`a un
$\bf
{C}$-bord par la m\'ethode de satellisation. \par

Soit $K$ un n\oe ud orient\'e dans $S^3$ et $J$ un n\oe ud orient\'e dans
un tore
solide standard, non nou\'e et plong\'e dans $S^3$. On peut aussi
consid\'erer $J$
comme un n\oe ud orient\'e dans $S^3$. Soit $N(K)$ un voisinage tubulaire
de $K$ et $f
: V \to N(K)$ un diff\'eomorphisme pr\'eservant l'orientation et qui envoie
la
longitude standard de $V$ sur celle du n\oe ud $K$. On note alors $K(J)$ le
n\oe ud
orient\'e dans $S^3$ obtenu comme l'image $f(J)$ de $J$ par $f$ dans $N(K)
\subset
S^3$. \par

On appelle ordre de $J$ dans $V$, le nombre minimal $t$ de points
d'intersection~de~$J$ avec un disque m\'eridien de $V$. On appelle nombre
de tours de
$J$ dans $V$,  le nombre $\omega$ d'intersection alg\'ebrique de $J$ avec
un
disque m\'eridien de $V$. On dit alors que $K(J)$ est obtenu par
satellisation de $J$
le long de $K$ avec ordre $t$ et nombre de tours $\omega$. \par

On obtient la condition n\'ecessaire suivante pour qu'un n\oe ud obtenu par
satellisation ne soit pas concordant \`a un $\bf {C}$-bord :
\bigskip

\noindent {\bf 2.5. Th\'eor\`eme.} {\it Soit $K$ un n\oe ud orient\'e dans
$S^3$. Soit
$J$ un n\oe ud orient\'e dans un tore solide non nou\'e standard $V$,
d'ordre $t$ et de
nombre de tours $\omega$. On suppose que $J^*$ en tant que n\oe ud dans
$S^3$ est
concordant \`a un $\bf {C}$-bord et que l'on a l'in\'egalit\'e $M(J) >
|\omega| M(K) +
{1 \over 2} (t - |\omega|)+1$. Alors le n\oe ud orient\'e $K(J)$, obtenu
par
satellisation de $J$ le long de $K$ avec ordre $t$ et nombre de tours
$\omega$, n'est
pas concordant \`a un $\bf {C}$-bord.}   \bigskip

{\bf Preuve.}~~Supposons que $K(J)$ est concordant \`a 
un $\bf {C}$-bord.~~Puisque $J^*$ est 
concordant \`a un $\bf {C}$-bord, la proposition 1.4 implique que
$M\bigl(K(J) \# J^*\bigr) \geq M\bigl(K(J)\bigr) + M(J^*)-1$. D'autre part,
on peut
effectuer la somme connexe avec $J^*$ dans une petite boule contenue dans
le
voisinage tubulaire $N(K)$ de $K$, d'o\`u $K(J)$ est \'equivalent \`a un
n\oe ud
obtenu en satelisant $J \# J^*$ autour de $K$ avec un ordre $t$ et un
nombre de
tours $\omega$, que l'on notera $K(J \# J^*)$. \par

D'apr\`es ${\rm [Sh],}$ puisque $J \# J^*$ est nul concordant, on a
l'in\'egalit\'e
$M\bigl(K(J \# J^*)\bigr) \break \leq |\omega| M(K) + {1 \over 2} (t -
|\omega|)$.
D'o\`u l'in\'egalit\'e $M\bigl(K(J)\bigr) + M(J^*)  \leq |\omega| M(K) + {1
\over 2}
(t - |\omega|)+1$, qui contredit l'hypoth\`ese puisque $M(J^*) = M(J)$ et
$M\bigl(K(J)\bigr) \geq 0$. \qquad \qquad \cqfd \par

Nous allons appliquer le Th\'eor\`eme 2.5 \`a des satellisations
particuli\`eres. \par

Rappelons que pour $p \in {\bf Z}$ et $q \geq 2$, le $(p, q)-$cable d'un
n\oe ud
orient\'e $K$ est le n\oe ud $K(p,q)$ obtenu par satellisation du n\oe ud
torique
$T(p, q)$ autour du n\oe ud $K$ avec ordre et nombre de tours $t = |\omega|
= q \geq
2$. De plus, $T^*(p,q) = T(-p, q)$. 
\bigskip

\noindent {\bf 2.6. Corollaire.} {\it Soient $K$ un n\oe ud orient\'e
quelconque dans
$S^3$, $p$ et $q$ deux entiers $\geq 2$  tels que $(p-1)(q-1) > 2qM(K)+1$.
Alors le
$(-p,q)-$cable, $K(-p,q)$,~de~$K$ n'est jamais concordant \`a un $\bf
{C}$-bord. En
particulier il existe des n\oe uds toriques it\'er\'es qui ne sont pas
concordants
\`a un $\bf {C}$-bord.
}\bigskip

{\bf Preuve.} D'apr\`es $[KM_{1,2}]$,  (cf. Th\'eorem. 1.3),
$M\bigl(T_{(p,q)}\bigr) =
{1 \over 2} (p-1)(q-1)$. Comme $t = |\omega| =q$, la condition $(p-1)(q-1)
>
2qM(K)+1$ \'equivaut \`a celle du Th\'eor\`eme 2.5. \qquad \qquad \cqfd
\bigskip

Nous allons maintenant consid\'erer le cas des n\oe uds doubles. Rappelons
que le
double (au sens de Whitehead [Wh]), $D(K \rho, \varepsilon)$, d'un n\oe ud
$K$, 
correspond au n\oe ud satellite $K(W_{\rho, \varepsilon})$, o\`u le n\oe ud
$W_{\rho, \varepsilon}$ dans le tore solide standard $V$ correspond au bord
d'un disque avec un clasp de signe $\varepsilon = +$ ou $-$ et admettant
$\rho$ tours
complets comme sur la figure.

{\BoxedEPSF{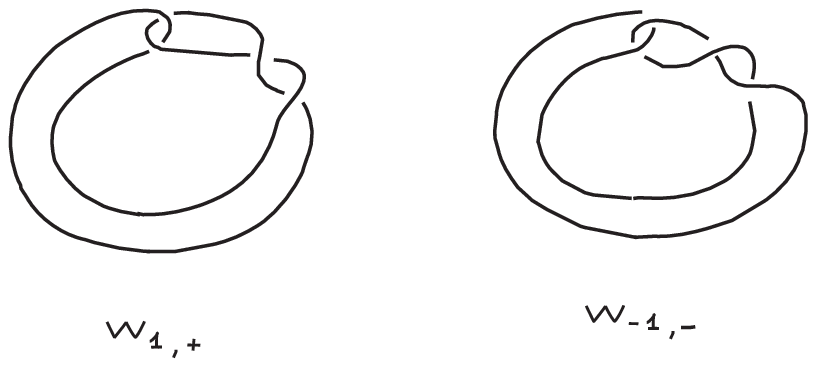}}

D'apr\`es $[R_8]$, le n\oe ud orient\'e $W_{\rho, \varepsilon}$ est
quasipositif,
donc un $\bf {C}$-bord, lorsque $\varepsilon = +$ et $\rho < 0$. De plus le
second
auteur a montr\'e $[R_5]$ que pour tout n\oe ud orient\'e $K$, il existe un
entier
$q(K)$ tel que pour tout entier $\rho \leq q(K)$, le n\oe ud orient\'e $D(K
; \rho ;
+)$ est quasipositif $\bigl(q(K)$ est le module de quasipositivit\'e de
$K$, ou bien
l'invariant de Thurston Bennequin $TB(K), cf. [R_{13}] \bigr)$. \par

Dans le cas des doubles n\'egatifs $D(K ; \rho, -)$ nous avons :
\bigskip

\noindent {\bf 2.7.~~Corollaire.}~~{\it Le double n\'egatif $D(K ; \rho, -)$
d'un n\oe ud
orient\'e nul concordant  $K$, avec $\rho > 0$ tours, n'est jamais
concordant
\`a un $\bf {C}$-bord.}
\medskip

{\bf Preuve.} Le n\oe ud oriente $W\rho, \ -\rho > 0$ est un n\oe ud qui a
pour
signature $\sigma(\omega_{\rho,-}) = -2$. Il n'est donc pas  nul
concordant. Comme $D(K ; \rho, -)$ et $W_{\rho,-}$ ont les m\^emes
invariants
d'Alexander, il en d\'ecoule que $D(K ; \rho, -)$ n'est pas  nul
concordant, et donc $M\bigl(D(K ; \rho,-)\bigr) = 1$. \par

Si on suppose que $D(K ; \rho, -)$ est concordant \`a un $\bf {C}$-bord
pour $\rho >
0$, la preuve du Th\'eor\`eme 2.4 montre que : $M\bigl(D(K ; \rho,-) \#
W^*_{\rho,-}\bigr) \geq M\bigl(D(K ; \rho, -)\bigr) + M(W^*_{\rho,-})-1$
car $W^*_{\rho,-} = W_{-\rho, +}$ est un $\bf{C}$-bord $[R_5]$. \par

Comme $W_{\rho, -}$ n'est pas nul concordant, $M(W_{\rho,-}) =
M(W^*_{\rho,-}) = 1$. D'o\`u \hfill\break $M\bigl(D(K ; \rho_{,-}) \#
W^*_{\rho,-}\bigr) \geq 1$. \par

D'un autre c\^ot\'e $D(K ; \rho_{,-}) \# W^*_{\rho,-}$ est obtenu par
satellisation
d'ordre $2$ et de nombre de tours $0$ du n\oe ud $W_{\rho,-} \#
W^*_{\rho,-}$  le long du n\oe ud nul concordant $K$.\par

Il d\'ecoule de $[Sh]$ que : $M\bigl(D(K ; \rho, -) \# W^*_{\rho,-}\bigr) = 0$, ce
qui donne la contradiction. \qquad \qquad \cqfd

\vskip .85truecm

\noindent \S\kern.15em 3 - {\bf Une obstruction \`a \^etre un
{\it spc}-${\bf C}$-bord.}   
\medskip

Pour l'instant, tous les entrelacs ${\bf C}$-bords connus v\'erifient en
fait une condition plus restrictive qui en fait des 
{\it spc}-${\bf C}$-bords.
\medskip

\noindent {\bf 3.1. D\'efinition.} {\it On dit qu'un entrelacs orient\'e
$L$ est un 
{\it spc}-${\bf C}$-bord s'il est obtenu comme intersection transverse d'une
courbe
alg\'ebrique plane dans ${\bf C}^2$ et du bord strictement pseudo-convexe
d'une boule
lisse $B^4 \subset {\bf C}^2$.}
\bigskip

Bien qu'on ne connaisse aucun exemple explicite, il serait \'etonnant que
tout
\hbox{${\bf C}$-bord} soit un {\it spc}-${\bf C}$-bord. On ne 
connait aucune obstruction
num\'erique pour 
\hbox{un entrelacs orient\'e} \`a \^etre un ${\bf C}$-bord. Par contre, la
condition de
stricte pseudo-convexit\'e du bord de la boule $B^4$ dans la construction
des
{\it spc}-${\bf C}$-bords permet d'\'etendre au cas des {\it spc}-${\bf C}$-bords une
obstruction
d\'ej\`a connue pour les entrelacs quasipositifs (cf. ${\rm [R_9]).}$ Ce
r\'esultat
conforte en particulier la conjecture suivante du second auteur :  \bigskip

\noindent {\bf Conjecture.} {\it Un entrelacs orient\'e qui est un
{\it spc}-${\bf C}$-bord
admet une pr\'esentation sous forme de tresse ferm\'ee quasipositive.} \par

Dans ce paragraphe on d\'emontre cette conjecture pour certains entrelacs
arborescents altern\'es, et en particulier pour les entrelacs \`a 2 ponts.
\par

L'obstruction \`a \^etre un {\it spc}-${\bf C}$-bord utilise le polyn\^ome \`a 2
variables
d'un entrelacs orient\'e (HOMFLY), nous commen\c cons par rappeler sa
d\'efinition en
utilisant pour les variables les notations de [Mo] (cf. ${\rm[R_9]
[R_{10}]).}$ \par

Soit $L$ un entrelacs orient\'e dans $S^3$. On note $P_L(v,z) \in {\bf
C}[v^{\pm1},
z^{\pm1}]$ le polyn\^ome de Laurent \`a 2 variables (ou polyn\^ome HOMFLY)
associ\'e
\`a l'entrelacs $L$ et d\'efini r\'ecursivement par : \par

{\BoxedEPSF{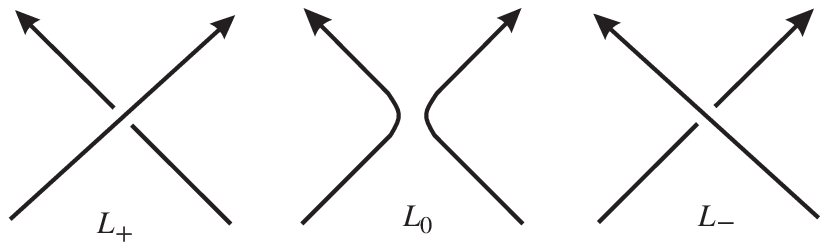}

\noindent $L_+, L_-$ et $L_0$ sont trois entrelacs dont les projections
sont
identiques, except\'e en un seul croisement comme indiqu\'e sur la Figure.
\par

$P_0(v, z)  =  1$, si  $0$ est \ le n\oe ud \ trivial \par
$P_{L_+} (v,z)  =  vz P_{L_0}(v,z) + v^2 P_{L_-}(v,z)  $. \par

En consid\'erant $P_L(v,z)$ comme un polyn\^ome de Laurent en $v$ \`a
coefficients
dans ${\bf C}[z^{\pm 1}]$, on obtient : $P_L(v, z) = \Sigma^D_{i=d} g_i(z)
v^i$,
$g_i(z) \in {\bf C} [z^{\pm 1}]$. \par

On d\'efinit alors $Ord_v P_L = d$, la valuation en $v$ de $P_L$. 
Le Th\'eor\`eme suivant 
permet de relier $Ord_v P_L$ et le grand genre de Murasugi $M(L)$,
lorsque
$L$ est un {\it spc}-${\bf C}$-bord.
\bigskip

\noindent {\bf 3.2. Th\'eor\`eme }[BR]{\bf.} {\it Soit $L$ un entrelacs orient\'e dans
$S^3$
ayant $r$ composantes connexes. Si $L$ est un {\it spc}-${\bf C}$-bord, alors
$Ord_v P_L
\geq 1-r + 2 M(L)$.}
\medskip

{\bf Preuve.} Soit $L = \partial B^4 \cap V$ un {\it spc}-${\bf C}$-bord, o\`u $V
\subset
{\bf C}^2$ est une courbe alg\'ebrique transverse au bord strictement
pseudo-convexe
de la boule $B^4 \subset {\bf C}^2$. \par

On peut toujour supposer, apre\`es une petite perturbations,
que $V$ est lisse.
Puisque $B^4$ est \`a bord strictement pseudo-convexe, la structure
complexe de ${\bf C}^2$ 
induit sur $\partial B^4$ la structure de contacte canonique [El]. Par
construction,
$L$ est un entrelacs transverse \`a cette structure de contact. On obtient
ainsi
(apr\`es une isotopie transverse [Be, Thm 10]) une pr\'esentation en tresse
ferm\'ee
$\widetilde \beta$ de $L$ qui\break v\'erifie 
la relation (cf. [La], [Be], [Sch]) :
$e(\widetilde \beta) - n(\widetilde \beta) = - \chi(V \cap B^4)$, o\`u
$e(\widetilde \beta)$ 
\hbox{est la longueur alg\'ebrique} de la tresse $\widetilde
\beta$ et
$n(\widetilde \beta)$ son nombre de brins. D'apr\`es\break%
Krohneimer-Mrowka ${\rm [KM_{1,2}]}$ (cf. Th\'eor\`eme
1.3)
$\chi(V \cap B^4)$ est la plus grande des caract\'eristiques
d'Euler
des surfaces lisses proprement plong\'ees dans $B^4$ et de bord $L$. D'o\`u, $\chi(V
\cap B^4) - \Sigma \mu_i = r - 2M(L)$. \par

D'apr\`es l'in\'egalit\'e de Frank-Williams [FW] et Morton [Mo], $Ord_v P_L
\geq
e(\widetilde \beta) - n(\widetilde \beta) + 1$, pour toute pr\'esentation
en tresse
ferm\'ee $\widetilde \beta$ de $L$. On en d\'eduit alors l'in\'egalit\'e
cherch\'ee.
\qquad \qquad \qquad   \cqfd 
\medskip

Le Corollaire suivant est imm\'ediat :
\medskip

\noindent {\bf 3.3.~~Corollaire.}~~{\it Soit $L$ un entrelacs orient\'ee dans
$S^3$
ayant $r$ composantes de bord. Si $L$ est un {\it spc}-${\bf C}$-bord, alors
$Ord_v P_L
\geq 1 - r$. } \qquad \qquad \qquad \cqfd
\bigskip

Dans le cas d'un n\oe ud $K$, si $\sigma(K)$ d\'esigne la signature de $K$,
d'apr\`es
K. Murasugi [Mu], on a $|\sigma(K)| \leq 2M(K)$. D'o\`u, on obtient :
\bigskip

\noindent {\bf 3.4. Corollaire.} {\it Soit $K$ un n\oe ud orient\'e dans
$S^3$. Si $K$ est un {\it spc}-${\bf C}$-bord, alors $Ord-v \ P_K \geq
|\sigma(K)|$. } \qquad
\qquad \qquad \cqfd

\bigskip
Nous appliquons maintenant la Proposition 3.2 \`a une famille d'entrelacs
orient\'es
arborescents altern\'es qui contient en particulier les entrelacs \`a 2
ponts. \par
 
\medskip

Rappelons d'abord l'op\'eration de plombage. 
Dans ce qui suit on note $A(0; n)$ la
surface de Seifert form\'ee par une bande orient\'ee d'\^ame non nou\'ee et
tordue par
$2n$ demi-tours positifs si $n \geq 0$ et n\'egatifs si $n \leq 0$. Le bord
de la
bande porte l'orientation induite par celle de la bande. Par exemple, $A(0,
-1)$
d\'esigne la bande de Hopf dont le bord est l'entrelacs alg\'ebrique de
Hopf. \par

Etant donn\'e une surface de Seifert $F \subset S^3$ orient\'ee et $\alpha
\subset F$
un arc proprement plong\'e (i.e. $\alpha \cap F = \partial \alpha)$, on
note $F * A(0
; n)$ la surface orient\'ee obtenue par plombage de la bande $A(0 ; n)$ sur
$F$ le
long de l'arc $\alpha$, par la construction suivante : \par

Soit $B^3 \subset S^3$ une boule telle que $F \subset B^3$ et $F \cap B^3 =
N(\alpha)$, un voisinage r\'egulier de $\alpha$ dans $F$. On recole alors
la bande
$A(0 ; n)$ \`a la surface $F$ le long de $N(\alpha)$ de telle sorte que :

\decale{\bf i)} $A(0 ; n) \subset S^3 - {\mathop B}^3$
\decale{\bf ii)} $A(0 ; n) \cap \partial B^3 = N(\alpha)$
\decale{\bf iii)} $\partial A(0 ; n) \cap \partial B^3 = \overline{\partial
N(\alpha)
- \partial F}$
\decale{\bf iv)} chaque composante connexe de $\partial F \cap N(\alpha)$
est un arc
essentiel dans $A(0 ; n)$. \par

Etant donn\'e un arbre planaire pond\'er\'e par des poids pairs, on peut
lui associer
une surface de Seifert orient\'ee bien d\'efinie, \`a isotopie pr\`es, dans
$S^3$ et
obtenue par plombage successifs de bandes tordues $A(0 ; n)$, $n \in {\bf
Z}$, comme
suit : \par

- \`a chaque sommet pond\'er\'e $e$ de poids $2n(e)$, on associe une bande
tordue
\break $A(0; n(e))$. \par
- \`a chaque ar\^ete correspond une op\'eration de plombage des deux bandes
tordues,
associ\'ees aux sommets de l'ar\^ete, le long d'un arc essentiel sur
chacune des
bandes (cf. Figure).

{\BoxedEPSF{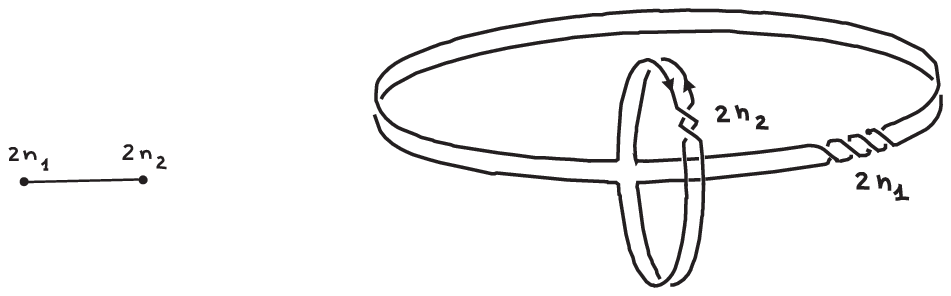}


Etant donn\'e un arbre planaire pond\'er\'e $T$, on consid\`ere la famille
$T_1, \cdots, T_k$\break%
des sous-arbres de $T$ obtenus en supprimant toutes les
ar\^etes de $T$
qui joignent deux sommets ayant des poids de signes oppos\'es. Tous les
sommets de
chacun des~arbres $T_i$, $i = 1, \cdots, k$, ont des poids de m\^eme signe.
On appelle
$\{T_1, \cdots, T_k\}$ une d\'ecomposition uniforme de $T$. On peut alors
toujours
supposer que les arbres $T_i, i = 1, \cdots, s$, ont tous des sommets de
poids
positifs, tandis que ceux des arbres $T_j, j = s + 1, \cdots, k$, sont
n\'egatifs. \par

On dit qu'un arbre planaire $T$, pond\'er\'e par des poids pairs, est {\bf
fortement
excessif} s'il existe une d\'ecomposition uniforme $\{T_1, \cdots, T_k\}$
de $T$ telle
que les poids $2n(e)$ de chaque arbre $T_i, i = 1, \cdots, k$, v\'erifient
: \par

a) $n(e) \not= 0$ pour chaque sommet $e \in T_i$ \par

b) $|n(e)| \geq v(e) - 1$ o\`u $v(e)$ est la valence du sommet $e$ dans
$T_i$. \par

Par exemple tout arbre lin\'eaire, pond\'er\'e par des poids pairs et non nuls, est
fortement excessif.
\bigskip

\noindent {\bf 3.5. D\'efinition.}  {\it On appelle entrelac arborescent
pair,
fortement excessif, l'entre- \break lacs orient\'e qui est le bord de la
surface de
.Seifert orient\'ee, obtenue par plombages successifs de bandes tordues,
suivant un
arbre planaire pond\'er\'e par des poids pairs et fortement excessif. (cf.
[MP, chap.
III]).} 
\medskip

Par exemple, tout entrelacs orient\'e \`a 2 ponts est un entrelacs
arborescent pair
fortement excessif. En particulier tout entrelacs arborescent pair et
fortement
ex\-cessif est altern\'e (cf. [MP, chap. III]). \par

Pour cette famille d'entrelacs orient\'es, on a une caract\'erisation
des {\it spc}-${\bf C}$-bord.
\bigskip

\noindent {\bf 3.6. Proposition.} {\it Un entrelacs arborescent pair et
fortement
excessif est un {\it spc}-${\bf C}$-bord si et seulement si tous les poids de
l'arbre
pond\'er\'e associ\'e sont strictement n\'egatifs.} 
\medskip

{\bf Preuve.} D'apr\`es ${\rm[R_8],}$ un entrelacs arborescent pair $L$ est
fortement
quasipositif si et seulement si  tous les poids de l'arbre pond\'er\'e
associ\'e sont
strictement n\'egatifs. En particulier un tel entrelacs orient\'e est un
{\it spc}-${\bf C}$-bord. \par

Supposons que $L$ est un entrelacs orient\'e arborescent pair et fortement
excessif,
associ\'e \`a un arbre pond\'er\'e $T$. D'apr\`es [MP, chap. III, 12,21] on
a : $Ord_v
P_L = p + q - 2 \displaystyle{\sum_{n_i>0}} n_i - 2s$, o\`u $p$ est le
nombre de
sommets de $T$ de poids n\'egatifs, $q$ celui des sommets de poids
positifs, $2n_i
> 0$ sont les poids positifs de $T$ et $s$ est le nombre de sous-arbres de la
d\'ecomposition uniforme de $T$ dont les sommets portent des poids
positifs. \par

D'autre part, puisque $L$ est un entrelacs arborescent pair, on a d'apr\`es
${\rm[R_8]}$ : $M(L) \geq {1 \over 2} (r-1 + p - q)$. \par

Si $L$ est un {\it spc}-${\bf C}$-bord, il d\'ecoule de la Proposition 3.2 que :
$Ord_v P_L
\geq p - q$. \par

On en d\'eduit alors l'in\'egalit\'e : $2q - 2 \displaystyle{\sum_{n_i>0}}
n_i - 2s \geq 0$. \par
\noindent Puisque l'arbre $T$ est fortement excessif, pour chaque sommet de
poids
positifs $2n_i \geq 1$, d'o\`u $2q - 2 \displaystyle{\sum_{n_i>0}} n_i \leq
0$. Il en
r\'esulte qu'en fait $s = 0$, et que tous les poids de l'arbre pond\'er\'e
$T$ sont
donc strictement n\'egatifs. \qquad \qquad \qquad \qquad \cqfd
\medskip

On obtient en particulier comme corollaire :
\medskip

\noindent {\bf 3.7.~~Corollaire.}~~{\it Un entrelacs arborescent pair
fortement excessif
(en particulier un entrelacs orient\'e \`a 2 ponts) est un 
{\it spc}-${\bf C}$-bord si et
seulement s'il est fortement quasipositif.}

\vskip .7truecm

\noindent {\bf BIBLIOGRAPHIE.}
\vskip 0,5 truecm
\item {[Be]} D. BENNEQUIN, ``Entrelacements et \'equations de Pfaff'',
Ast\'erisque 107-108, (1983), 87-161.

\item {[BR]} M. BOILEAU et L. RUDOLPH, ``Une obstruction topologique
\`a ce qu'un n\oe ud dans le bord strictement pseudo-convexe d'une
boule $B^4$ borde un morceau de courbe alg\'ebrique'', \`a para\^\i tre
(1995).

\item {[BW]} M. BOILEAU et C. WEBER, ``Le Probl\`eme de J. Milnor sur le
nombre
gordien des n\oe uds alg\'ebriques'', 
Mono. n$^{\scriptsize\circ}$ 31 de l'Enseignement
Math., N\oe uds, Tresses et Singularit\'es, (1983), 49-98.

\item {[CG]} T.D. COCHRAN et R.E. GOMPF, ``Application of Donaldson's
theorems to classical knot concordance, homology 3-spheres 
and Property P'', Topology 27,
(1988), 495-512.

\item {[EN]} D. EISENBUD et W. NEUMANN, ``Three-dimensional link theory and
invariants of plane curve singularities'', Annals of Math. Studies 110,
(1985),
Princeton University Press.

\item {[El]} Ya. ELIASHBERG, ``Filling by holomorphic disks and its
applications'', Lon-\break%
don Math. Soc. Lect. Notes Ser. 151 (1991),
45-67.

\item {[F]} M. FREEDMAN, ``A Surgery sequence in dimension 4 ; the relations
with knot concordance'', Invent. Math. 68 (1982), 195-226.

\item {[FW]} J. FRANKS et R. WILLIAMS, ``Braids and the Jones-Conway
polynomial'', T.A.M.S. 303 (1987), 97-108.

\item {[Go]} R.E. GOMPF, ``Smooth concordance of topologically slice knots'',
Topology 25 (1986), 353-373.

\item {[KL]} R. KIRBY et W. LICKORISH ``Prime knots and concordance'' Math.
Proc.
Cambridge Philos. Soc. 86 (1979), 437-441.

\item {[KM$_{1,2}$]} P. KRONHEIMER et T. MROWKA, ``Gauge theory for embedded
surfaces'',~{\bf I}\break%
: Topology 32 (1993), 773-826 -  {\bf II} : to appear in
Topology. 

\item {[KM$_3$]} P. KRONHEIMER et T. MROWKA, 
``The genus of embedded surfaces in
the projective plane'', Math. Res. Letters 1 (1994), 797-808.

\item {[La]} H. LAUFER, ``On the number of singularities of an analytic
curve'', T.A.M.S. 186 (1969), 527-535.

\item {[L\^e]} LE D.T., ``Sur les n\oe uds alg\'ebriques'', Compositio Math.
25
(1972), 281-321.

\item {[Mi]} J. MILNOR, ``Singular points of complex hypersurfaces'', Ann. of
Math. Stud-\break%
ies 61 (1968), Princeton University Press.

\item {[Mo]} H. MORTON ``Seifert circles and knot polynomials'', Math. Proc.
Camb. Phil. Soc. 99 (1986), 107-110.

\item {[Mu]} K. MURASUGI ``On certain numerical invariant of link types''
T.A.M.S. 117 (1965), 387-422.

\item {[MP]} ``K. MURASUGI et J.H. PRZYTYCKI, ``An Index of a Graph with
Applications to Knot Theory'' Memoirs of A.M.S. 508 - Vol. 106 (1993).

\item {[Ne]} W. NEUMANN, ``Complex algebraic plane curves via their links at
infinity'', Invent. Math. (1989), 445-489.

\item {[NR]} W. NEUMANN et L. RUDOLPH, ``Unfoldings in knot theory'', Math.
Annalen 278 (1987), 409-439 ; corrigendum ibid. (1988).

\item {[R$_1$]} L. RUDOLPH, ``Embeddings of the line in the plane'', J. Reine
Ange. Math. 337 (1982), 113-118.

\item {[R$_2$]} L. RUDOLPH, ``Algebraic functions and closed braids'',
Topology
22 (1983), 191-202.

\baselineskip=14pt plus 1pt minus 1pt
\item {[R$_3$]} L. RUDOLPH, ``Braided surfaces and Seifert ribbons for
closed
braids'', Comm. Math. Helv. 58 (1983), 1-37.

\item {[R$_4$]} L. RUDOLPH, ``Construction of quasipositive knots and links
I'', Mono. n$^{\scriptsize\circ}$\break%
31 de l'Enseignement
Math., N\oe uds, Tresses et Singularit\'es, (1983), 
233-345.

\item {[R$_5$]} L. RUDOLPH, ``Construction of quasipositive knots and links
II,
Contemp. Math. 35 (1984), 485-491.

\item {[R$_6$]} L. RUDOLPH, ``A characterisation of quasipositive Seifert
surfaces (Constuctions of quasipositive knots and links III), Topology 31
(1992), 231-237.

\item {[R$_7$]} L. RUDOLPH, ``Quasipositive Annuli (Constructions of
quasipositive knots\break%
and links IV), J. Knot Theory Ramif 1 (1993), 451-466.

\item {[R$_8$]} L. RUDOLPH, ``Quasipositive plumbing (Constructions of
quasipositive knots and links V), preprint (1990, 1995).

\item {[R$_9$]} L. RUDOLPH, ``A congruence between link polynomials'', Math.
Proc. Camb. Phil. Soc. 107 (1990), 319-327.

\item {[R$_{10}$]} L. RUDOLPH, ``Quasipositivity and new knot invariants''
Rev.
Mat. Univ. Complutense Madrid 2 (1989), 85-109.

\item {[R$_{11}$]} L. RUDOLPH, ``Special positions for surfaces bounded by
closed braids'', Rev. Math. Iberoamericana 1 (1985), 93-133.

\item {[R$_{12}$]} L. RUDOLPH, ``Quasipositivity as an obstruction to
sliceness'', Bull. Am.\break%
Math. Soc. 29 (1993), 51-59.

\item {[R$_{13}$]} L. RUDOLPH, ``An obstruction to sliceness via contact
geometry and `classical' gauge theory'', Invent. Math. 119 (1995), 155-163.

\item {[R$_{14}$]} L. RUDOLPH, ``Totally tangential links of intersection of
complex plane curves with round spheres'', Topology '90, de Gruyter (1992),
343-349.

\item {[Sch]} N.V. SCHERBINA, ``Traces of pluriharmonic functions on the
boundaries of analytic varieties'', Math. Z. 213 (1993), 171-177.

\item {[Shi]} T. SHIBUYA, ``On 4-dimensional genera of compound knots'', 
Kobe Math.\break%
Sem. Notes 8 (1980), 299-305.

\item {[St]} G.P.STOLZENBERG, ``Uniform approximation on smooth curves'',
Acta Math. 115 (1966), 185-198.

\item {[Wh]} J.H.C. WHITEHEAD, ``On doubled knots'', J. London Math. Soc. 12
(1937), 63-71.
\eject

M. BOILEAU, Laboratoire Topologie et G\'eometrie, URA CNRS 1408, 
Universit\'e Paul Sabatier 31062 TOULOUSE C\'edex France

L. RUDOLPH, Department of Mathematics and Computer Science, Clark 
University, 950 Main Street, Worcester, Massachusetts 01610 USA
\vfill
\eject
\end